\documentstyle[11pt,epic]{article}
\begin{document}
\newcommand{\up}{\vspace*{-0.2cm}}
\newcommand{\upp}{\vspace*{-0.3cm}}
\newcommand{\qed}{\hfill$\rule{.05in}{.1in}$\vspace{.3cm}}
\newcommand{\pf}{\noindent{\bf Proof: }}
\newtheorem{thm}{Theorem}
\newtheorem{lem}{Lemma}
\newtheorem{prop}{Proposition}
\newtheorem{prob}{Problem}
\newtheorem{ex}{Example}
\newtheorem{cor}{Corollary}
\newtheorem{conj}{Conjecture}
\newtheorem{cl}{Claim}
\newtheorem{df}{Definition}
\newtheorem{rem}{Remark}
\newcommand{\beq}{\begin{equation}}
\newcommand{\eeq}{\end{equation}}
\newcommand{\<}[1]{\left\langle{#1}\right\rangle}
\newcommand{\be}{\begin{enumerate}}
\newcommand{\ee}{\end{enumerate}}
\newcommand{\Bul}{\mbox{$\bullet$ } }
\newcommand{\al}{\alpha}
\newcommand{\ep}{\epsilon}
\newcommand{\si}{\sigma}
\newcommand{\om}{\omega}
\newcommand{\la}{\lambda}
\newcommand{\La}{\Lambda}
\newcommand{\Ga}{\Gamma}
\newcommand{\ga}{\gamma}
\newcommand{\im}{\Rightarrow}
\newcommand{\2}{\vspace{.2cm}}
\newcommand{\es}{\emptyset}

\title{\bf Discrepancy and Signed Domination in Graphs and Hypergraphs}
\author{A. Poghosyan and V. Zverovich\\
\\
{\footnotesize Department of Mathematics and Statistics}\\
{\footnotesize University of the West of England}\\
{\footnotesize Bristol, BS16 1QY}\\
{\footnotesize UK}}
\date{}
\maketitle

\begin{abstract}
For a graph $G$, a signed domination function of $G$ is a
two-colouring of the vertices of $G$ with colours +1 and --1 such
that the closed neighbourhood of every vertex contains more +1's
than --1's. This concept is closely related to combinatorial
discrepancy theory as shown by F\"{u}redi and Mubayi
[{\em J. Combin. Theory}, Ser. B {\bf 76} (1999) 223--239]. The
signed domination number of $G$ is the minimum of the sum of
colours for all vertices, taken over all signed domination
functions of $G$. In this paper, we present new upper and lower
bounds for the signed domination number. These new bounds improve
a number of known results.

\vspace*{.2cm} \noindent {\footnotesize Keywords: {\it graphs,
signed domination function, signed domination number.}}
\end{abstract}

\section{Discrepancy Theory and Signed Domination} 

Originated from number theory, discrepancy theory is, generally
speaking, the study of irregularities of distributions in various
settings. Classical combinatorial discrepancy theory is devoted to
the problem of partitioning the vertex set of a hypergraph into
two classes in such a way that all hyperedges are split into
approximately equal parts by the classes, i.e. we are interested
in measuring the deviation of an optimal partition from perfect,
when all hyperedges are split into equal parts. It may be pointed
out that many classical results in various areas of mathematics,
e.g. geometry and number theory, can be formulated in these terms.
Combinatorial discrepancy theory was introduced and studied by
Beck in \cite{beck1}. Also, studies on discrepancy theory have
been conducted in \cite{beck3,beck2, beck4} and \cite{sos}.

Let ${\cal H}=(V,{\cal E})$ be a hypergraph with the vertex set $V$
and the hyperedge set ${\cal E}=\{E_1,...,E_m\}.$ One of the main
problems in classical combinatorial discrepancy theory is to colour
the elements of $V$ by two colours in such a way that all of the
hyperedges have almost the same number of elements of each colour.
Such a partition of $V$ into two classes can be represented by a
function
$$f:V\rightarrow \{+1, -1\}.$$
 For a set $E\subseteq V$, let us define the {\em imbalance} of $E$ as follows:
$$
f(E)=\sum_{v \in E} f(v).
$$
 First defined by Beck \cite{beck1}, the
{\em discrepancy of ${\cal H}$ with respect to} $f$ is
$$
{\cal D}({\cal H},f)=\max_{E_i \in {\cal E}}|f(E_i)|
$$
and the {\em discrepancy} of ${\cal H}$ is
$$
{\cal D}({\cal H})=\min_{\tiny f:\; V\rightarrow\{+1,-1\}}{\cal
D}({\cal H},f).
$$
Thus, the discrepancy of a hypergraph tells us how well all its
hyperedges can be partitioned.  Spencer \cite{spe1} proved a
fundamental ``six-standard-deviation" result that for any hypergraph
${\cal H}$ with $n$ vertices and $n$ hyperedges,
$$
{\cal D}({\cal H}) \le 6 \sqrt n.
$$
As shown in \cite{alo}, this bound is best possible up to a constant
factor. More precisely, if a Hadamard matrix of order $n>1$ exists,
then there is a hypergraph ${\cal H}$ with $n$ vertices and $n$
hyperedges such that
$$
{\cal D}({\cal H}) \ge 0.5 \sqrt n.
$$
It is well known that a Hadamard matrix of order between $n$ and
$(1-\epsilon)n$ does exist for any $\epsilon$ and sufficiently large
$n$. The following important result, due to Beck and Fiala
\cite{beck2}, is valid for a hypergraph with any number of
hyperedges:
$$
{\cal D}({\cal H}) \le 2\Delta -1,
$$
where $\Delta$ is the maximum degree of vertices of ${\cal H}$. They
also posed the discrepancy conjecture that for some constant $K$
$$
{\cal D}({\cal H}) < K \sqrt\Delta.
$$

Another interesting aspect of discrepancy was discussed by
F\"{u}redi and Mubayi in their fundamental paper \cite{fur}. A
function $g: V\rightarrow \{+1, -1\}$ is called a {\em signed
domination function} (SDF) of the hypergraph ${\cal H}$ if
$$
g(E_i)=\sum_{v \in E_i} g(v) \ge 1
$$
for every hyperedge $E_i\in \cal E$, i.e. each hyperedge has a
positive imbalance. The {\em signed discrepancy} of ${\cal H}$,
denoted by ${\cal SD}({\cal H})$, is defined in the following way:
$$
{\cal SD}({\cal H})=\min_{{\mathrm SDF} g} g(V),
$$
where the minimum is taken over all signed domination functions of
${\cal H}$. Thus, in this version of discrepancy, the success is
measured by minimizing the imbalance of the vertex set $V$, while
keeping the imbalance of every hyperedge $E_i\in \cal E$ positive.

One  of the main results in this context, formulated in terms of
hypergraphs, is due to F\"{u}redi and Mubayi \cite{fur}:

\begin{thm}[\cite{fur}] \label{f_m}
Let ${\cal H}$ be an $n$-vertex hypergraph with hyperedge set
${\cal E}= \{E_1,...,E_m\}$, and suppose that every hyperedge has
at least $k$ vertices, where $k \ge 100.$ Then
$$
{\cal SD}({\cal H}) \le 4{\sqrt{\ln k \over k}}n + {1 \over k}m.
$$
\end{thm}

This theorem can be easily re-formulated in terms of graphs by
considering the neighbourhood hypergraph of a given graph. Note that
$\ga_s(G)$ is the signed domination number of a graph $G$, which is
formally defined in the next section.

\begin{thm}[\cite{fur}] \label{fur_m_thm}
If $G$ has $n$ vertices and minimum degree $\delta \ge 99$, then
$$
\ga_s(G) \le \left(4\sqrt{{\ln(\delta+1) \over \delta+1}} + {1 \over \delta+1} \right)n.
$$
\end{thm}

Moreover, F\"{u}redi and Mubayi \cite{fur} found quite good upper
bounds for very small values of $\delta$ and, using Hadamard
matrices, constructed a $\delta$-regular graph $G$ of order
$4\delta$ with
$$
\ga_s(G) \ge 0.5 \sqrt \delta - O(1).
$$
This construction shows that the upper bound in Theorem
\ref{fur_m_thm} is off from optimal by at most the factor of
$\sqrt{\ln\delta}$. They posed an interesting conjecture that, for
some constant $C$,
$$
\ga_s(G) \le {C \over \sqrt\delta} n,
$$
and proved that the above discrepancy conjecture, if true, would
imply this upper bound for $\delta$-regular graphs. A strong result
of Matou\v{s}ek \cite{mat} shows that the bound is true, but the
constant $C$ in his proof is big making the result of rather
theoretical interest.


The lower bound for the signed domination number given in the
theorem below is formulated in terms of the degree sequence of a
graph. Other lower bounds are also known, see Corollaries
\ref{haas-zhang}, \ref{lb-Cor1} and \ref{lb-Cor2}.

\begin{thm}[\cite{dun1}] \label{dunb_low}
Let $G$ be a graph with degrees $d_1 \le d_2 \le...\le d_n.$ If $k$ is the smallest integer for which
$$
\sum_{i=0}^{k-1}d_{n-i} \ge 2(n-k) + \sum_{i=1}^{n-k}d_{i},
$$
then
$$\ga_s(G) \ge 2k-n.$$
\end{thm}

In this paper, we present new upper and lower bounds for the
signed domination number, which improve the above theorems and
also generalise three known results formulated in Corollaries
\ref{haas-zhang}, \ref{lb-Cor1} and \ref{lb-Cor2}. Note that our
results can be easily re-formulated in terms of hypergraphs.
Moreover, we rectify F\"{u}redi--Mubayi's conjecture formulated
above as follows: for some $C\le 10$ and $\alpha$, $0.18\le \alpha
< 0.21$,
$$
\ga_s(G) \le \min\Big\{{n\over \delta^\alpha}, {Cn\over
\sqrt\delta}\Big\}.
$$

\section{Notation and Technical Results}

All graphs will be finite and undirected without loops and multiple
edges. If $G$ is a graph of order $n$, then
$V(G)=\{v_1,v_2,...,v_n\}$ is the set of vertices in $G$ and $d_i$
denotes the degree of $v_i$. Let $N(x)$ denote the neighbourhood of
a vertex $x$. Also, let $ N(X)=\cup_{x\in X} {N(x)} $ and $
N[X]=N(X)\cup X$. Denote by $\delta(G)$ and $\Delta(G)$ the minimum
and maximum degrees of vertices of $G$, respectively. Put
$\delta=\delta(G)$ and $\Delta=\Delta(G)$.

A set $X$ is called a {\it dominating set} if every vertex not in
$X$ is adjacent to a vertex in $X$. The minimum cardinality of a
dominating set of $G$ is called the \emph{domination number}
$\gamma(G)$. The domination number can be defined equivalently by
means of a {\it domination function}, which can be considered as a
characteristic function of a dominating set in $G$. A function $f :
V(G) \rightarrow \{0, 1\}$ is a domination function on a graph $G$
if for each vertex ${v \in V(G)}$,
\begin{equation}\label{exfirst}
\sum_{x \in N[v]}f(x) \ge 1.
\end{equation}
 The value $\sum_{v \in V(G)} f(v)$
is called the weight $f(V(G))$ of the function $f$. It is obvious
that the minimum of weights, taken over all domination functions on
$G$, is the domination number $\gamma(G)$ of $G$.

It is easy to obtain different variations of the domination number
by replacing the set \{0,1\} by another set of numbers. If \{0,1\}
is exchanged by \{-1,1\}, then we obtain {\it the signed
domination number}. A signed domination function of a graph $G$
was defined in \cite{dun1} as a function $f : V(G) \rightarrow
\{-1, 1 \}$ such that for each $v \in V(G)$, the expression
(\ref{exfirst}) is true. The signed domination number of a graph
$G$, denoted $\gamma_s(G)$, is the minimum of weights $f(V(G))$,
taken over all signed domination functions $f$ on $G$. A research
on signed domination has been carried out in
\cite{dun1}--\cite{hen2} and \cite{mat}.

Let $d \ge 2$ be an integer and $0\le p \le 1.$ Let us denote
$$
f(d,p)=\sum_{m=0}^{\lceil 0.5d \rceil} (\lceil 0.5d \rceil-m+1)\pmatrix{d+1 \cr m}p^m(1-p)^{d+1-m}.
$$

We will need the following technical results:

\begin{lem}[\cite{fur}] \label{fur_m}
If $d$ is odd, then $$f(d+1,p) < 2(1-p)f(d,p).$$ If $d$ is even, then
$$f(d+1,p)<\left(2p+(1-p){d+4\over d+2}\right)f(d,p).$$
In particular, if $$2(1-p)\left(2p+(1-p){d+4 \over
d+2}\right)<1,$$ then \begin{displaymath}\max_{d \ge
\delta}{f(d,p)\in \{f(\delta,p),f(\delta+1,p)\}}.
\end{displaymath}
\end{lem}

\begin{lem}[\cite{cher}] \label{chernoff}
Let  $p \in [0,1]$ and $X_1,...,X_k$ be mutually independent random
variables with
\begin{eqnarray*}
&{\mathbf P}[X_i=1-p]=p, \\
&{\mathbf P}[X_i=-p]= 1-p.
\end{eqnarray*}
If $X= X_1+...+X_k$ and $c>0$, then
$$
{\mathbf P}[X < -c] < e^{-{c^{2} \over 2pk}}.
$$
\end{lem}

Let us also denote
$$
{\widetilde d}_{0.5} = \pmatrix{\delta'+1 \cr {\lceil 0.5 \delta'
\rceil}},
$$
where
\begin{displaymath}
\delta'= \cases{\\
\delta& if $\delta$ is odd;\cr\\
\delta+1& if $\delta$ is even.\\
}\\
\end{displaymath}

\section{Upper Bounds for the Signed Domination Number}

The following theorem provides an upper bound for the signed
domination number, which is better than the bound of Theorem
\ref{fur_m_thm} for `relatively small' values of $\delta$. For
example, if $\delta(G)=99$, then, by Theorem \ref{fur_m_thm},
$\ga_{s}(G) \le  0.869n,$ while Theorem \ref{main_signed} yields
$\ga_{s}(G) \le  0.537n$. For larger values of $\delta$, the
latter result is improved in Corollaries \ref{c1}--\ref{c3}.

\begin{thm}\label{main_signed}
 For any graph $G$ with $\delta > 1,$
\begin{equation}\label{boundsigned}
\ga_{s}(G) \le \left(1-{{2\widehat{\delta}} \over
{(1+{\widehat{\delta}})}^{1+1/{\widehat{\delta}}}\,{\widetilde
d}_{0.5}^{\;1/{\widehat{\delta}}}} \right) n,
\end{equation}
where
$ \widehat{\delta} = \lfloor 0.5 \delta \rfloor$.

\end{thm}

\pf Let $A$ be a set formed by an independent choice of vertices
of $G$, where each vertex is selected with the probability
$$p = 1- {1 \over (1+\widehat{\delta})^{1/\widehat{\delta}}\;{\widetilde
d}_{0.5}^{\;1/\widehat{\delta}}}.$$ For $m\ge0$, let us denote by $B_m$ the set of
vertices $v\in V(G)$ dominated by exactly $m$ vertices of $A$ and
such that $|N[v]\cap A|< \lceil 0.5 d_v \rceil +1$, i.e.
$$
|N[v]\cap A| = m \le \lceil 0.5 d_v \rceil.
$$
Note that each vertex $v\in V(G)$ is in at most one of the sets
$B_m$ and $0\le m\le \lceil 0.5 d_v \rceil$. Then we form a set $B$
by selecting $\lceil 0.5 d_v \rceil - m + 1 $ vertices from $N[v]$
that are not in $A$ for each vertex $v \in B_m$ and adding them to
$B$. We construct the set $D$ as follows: $D=A\cup B$. Let us assume
that $f$ is a function $f : V(G) \rightarrow \{-1, 1\}$ such that
all vertices in $D$ are labelled by $1$ and all other vertices by
$-1$. It is obvious that $f(V(G))=|D| - (n - |D|)$ and $f$ is a
signed domination function.

 The expectation of $f(V(G))$ is
\begin{eqnarray*}
{\mathbf E}[f(V(G))] &=& 2{\mathbf E}[|D|]- n\\
&=& 2({\mathbf E}[|A|] + {\mathbf E}[|B|]) - n\\
&\le& 2 \, \sum_{i=1}^{n}{\mathbf P}(v_i\in A)+ 2
\sum_{i=1}^{n}\sum_{m=0}^{\lceil
0.5 d_i \rceil}(\lceil 0.5 d_i \rceil - m + 1) {\mathbf P}(v_i\in B_m) - n \\
&=& 2pn + 2\sum_{i=1}^{n}\sum_{m=0}^{\lceil 0.5 d_i \rceil}(\lceil
0.5 d_i \rceil - m + 1) \pmatrix{d_i+1 \cr m} p^m (1-p)^{d_i+1-m}
- n\\
&\le& 2pn + 2\sum_{i=1}^{n}{\max_{d_i \ge \delta}f(d_i,p)-n}.
\end{eqnarray*}
It is not difficult to check that
$2(1-p)(2p+(1-p)(d+4)/(d+2)) < 1$ for any $d \ge \delta \ge 2$. By Lemma \ref{fur_m},
$$
\max_{d \ge \delta}{f(d,p) \in \{f(\delta,p),f(\delta+1,p)}\}.
$$
The last inequality implies $2(1-p)<1$.
Therefore, by Lemma \ref{fur_m},
$$
\max_{d \ge \delta}{f(d,p)}=f(\delta,p)
$$
if $\delta$ is odd. If $\delta$ is even, then we can prove that
$$
\max_{d \ge \delta}{f(d,p)}=f(\delta+1,p).
$$
Thus,
$$
\max_{d \ge \delta}{f(d,p)}=f(\delta',p).
$$
Therefore,
$$
{\mathbf E}[f(V(G))] \le 2pn + 2n\sum_{m=0}^{\lceil 0.5 \delta'
\rceil}(\lceil 0.5 \delta' \rceil - m + 1) \pmatrix{\delta'+1 \cr
m} p^m (1-p)^{\delta'+1-m} - n.
$$
Since
$$(\lceil 0.5 \delta' \rceil - m + 1) \pmatrix{\delta'+1 \cr m} \le \pmatrix{\delta'+1 \cr \lceil 0.5 \delta' \rceil} \pmatrix{\lceil 0.5 \delta' \rceil \cr m},$$
we obtain
\begin{eqnarray*}
{\mathbf E}[f(V(G))]&\le& 2pn + 2n\sum_{m=0}^{\lceil 0.5 \delta' \rceil} \pmatrix{\delta'+1 \cr \lceil 0.5 \delta' \rceil} \pmatrix{\lceil 0.5 \delta' \rceil \cr m} p^m (1-p)^{\delta'+1-m} - n\\
&=& 2pn + 2n\pmatrix{\delta'+1 \cr \lceil 0.5 \delta' \rceil} (1-p)^{\delta'-\lceil 0.5 \delta' \rceil+1} \sum_{m=0}^{\lceil 0.5 \delta' \rceil} \pmatrix{\lceil 0.5 \delta' \rceil \cr m} p^m (1-p)^{\lceil 0.5 \delta' \rceil-m} - n\\
&=& 2pn + 2n{\widetilde d}_{0.5}(1-p)^{\delta'-\lceil 0.5 \delta' \rceil+1} - n.
\end{eqnarray*}
Taking into account that $\delta'-\lceil0.5 \delta'\rceil=\lfloor0.5\delta'\rfloor=\lfloor0.5\delta\rfloor=\widehat{\delta},$ we have
\begin{eqnarray*}
{\mathbf E}[f(V(G))] &\le& 2pn + 2n{\widetilde
d}_{0.5}(1-p)^{\widehat{\delta}+1} - n\\
&\le& \left(1-{2\;\widehat{\delta} \over
{(1+\widehat{\delta})}^{1+1/\widehat{\delta}} \,{\widetilde
d}_{0.5}^{\;1/\widehat{\delta}}} \right) n,
\end{eqnarray*}
as required. The proof of Theorem \ref{main_signed} is complete.
\qed

Our next result and its corollaries give a modest improvement of
Theorem \ref{fur_m_thm}. More precisely, the upper bound of Theorem
\ref{main_2} is asymptotically 1.63 times better than the bound of
Theorem \ref{fur_m_thm}, and for $\delta=10^6$ the improvement is
1.44 times.

\begin{thm}\label{main_2}
 If $\delta(G) \ge 10^6$, then
$$
\ga_{s}(G) \le {\sqrt{6\ln(\delta+1)}+ 1.21 \over \sqrt{\delta +1}}n.
$$
\end{thm}

\pf Denote $\delta^{+}=\delta+1$, $N_v=N[v]$ and $n_v=|N_v|.$ Let $A$ be a set formed by an independent choice of vertices
of $G$, where each vertex is selected with the probability
$$p= 0.5 + \sqrt{1.5\ln{\delta^{+}/\delta^{+}}}. $$
Let us construct two sets $Q$ and $U$ in the following way: for each
vertex $ v \in V(G)$, if $|N_v\cap A| \le 0.5 n_v,$ then we put $v
\in U$ and add $\lfloor 0.5 n_v +1 \rfloor $ vertices of $N_v$ to
$Q$. Furthermore, we assign ``+" to $A\cup Q$, and ``--" to all
other vertices. The resulting function $g:V(G)\rightarrow \{-1, 1\}$
is a signed domination function, and
$$
g(V(G))=2|A \cup Q| - n \le 2|A| + 2|Q|-n.
$$
The expectation of $g(V(G))$ is
\begin{eqnarray}
\nonumber
{\mathbf E}[g(V(G))] &\le& 2{\mathbf E}[|A|]+2{\mathbf E}[|Q|]-n\\
\label{1}&=& 2pn - n + 2{\mathbf E}[|Q|].
\end{eqnarray}
It is easy to see that $|Q| \le \sum_{v \in \; U}\lfloor 0.5n_v+1\rfloor,$ hence
\begin{equation} \label{4}
{\mathbf E}[|Q|]\le\sum_{v \in V(G)}\lfloor
0.5n_v+1\rfloor\;{\mathbf P}[v \in U],
\end{equation}
where
$$
{\mathbf P}[v\in U]={\mathbf P}[|N_v \cap A| \le 0.5 n_v].
$$

Let us define the following random variables
\begin{displaymath}
X_w= \cases{\\
1-p& if $w \in A$\cr\\
-p& if $w \notin A$\\
}\\
\end{displaymath}
and let $X_v^{*}=\sum_{w \in N_v}X_w.$
We have
$$
|N_v \cap A| \le 0.5 n_v \mbox{\quad if and only if
\quad}X_v^{*}\le(1-p)0.5 n_v+(-p)0.5 n_v.
$$
Thus,
$$
{\mathbf P}[|N_v \cap A| \le 0.5 n_v]= {\mathbf P}[
X_v^{*}\le(0.5-p)n_v].
$$
By Lemma \ref{chernoff},
$$
{\mathbf P}[X_v^{*} \le \left(0.5 - p\right)n_v] <
e^{-{1.5n_v\ln{\delta^{+}/\delta^{+}} \over 1 +
\sqrt{6\ln{\delta^{+}/\delta^{+}}}}}.
$$
For $n_v \ge \delta^{+}>10^{6}$, let us define
$$
y(n_v, \delta^{+})={1.5n_v\ln{\delta^{+}/\delta^{+}} \over 1 + \sqrt{6\ln{\delta^{+}/\delta^{+}}}}-\ln(2.25n_v^{1.5})+1.
$$
The function $y(n_v, \delta^{+})$ is an increasing function of $n_v$ and
$y(\delta^{+},\delta^{+})>0$ for $\delta^+>10^6.$
Hence $y(n_v, \delta^{+})\ge y(\delta^{+}, \delta^{+})>0$ and
$$
{1.5n_v\ln{\delta^{+}/\delta^{+}} \over 1 +
\sqrt{6\ln{\delta^{+}/\delta^{+}}}} > \ln(2.25n_v^{1.5}) -1.
$$
We obtain
$${\mathbf P} [|N_v\cap A| \le 0.5 n_v] < e^{1-\ln(2.25n_v^{1.5})} = {e \over 2.25n_v^{1.5}},$$
and, using inequality (\ref{4}),
$$
2{\mathbf E}[|Q|]\le 2\sum_{v\in V(G)}{e(0.5n_v+1) \over
2.25n_v^{1.5}} \le {e(\delta +3)n \over 2.25(\delta+1)^{1.5}}\le
{1.21 \over \sqrt{\delta+1}}n.
$$
Thus, (\ref{1}) yields
\begin{eqnarray*}
{\mathbf E}[g(V(G))]&\le& 2pn - n + {1.21n \over \sqrt{\delta +1}}\\
&=& {\sqrt{6\ln(\delta+1)}+ 1.21 \over \sqrt{\delta +1}}n,
\end{eqnarray*}
as required. The proof of Theorem \ref{main_2} is complete. \qed \\

\begin{cor} \label{c1}
If $24,000 \le \delta,$ then
$$
\gamma_s(G) \le {\sqrt{6.8\ln(\delta+1)}+0.32 \over \sqrt{\delta+1}}n.
$$
\end{cor}

\pf Putting $p=0.5 + \sqrt{1.7\ln{\delta^+}/\delta^+}$ in the
proof of Theorem \ref{main_2}, we obtain by Lemma \ref{chernoff},
$$
{\mathbf P}[X_v^{*} \le \left(0.5 - p\right)n_v] <
e^{-{1.7n_v\ln{\delta^{+}/\delta^{+}} \over 1 +
\sqrt{6.8\ln{\delta^{+}/\delta^{+}}}}}.
$$
Let us define the following function:
$$
y(n_v, \delta^+)= {1.7n_v\ln\delta^+/\delta^+ \over 1+
\sqrt{6.8\ln\delta^+/\delta^+}} - \ln(3.14n_v^{1.5})
$$
for $n_v\ge \delta^+ > 24,000$. The function $y(n_v, \delta^{+})$ is
an increasing function of $n_v$ and $y(\delta^{+},\delta^{+})>0$ for
$\delta^+>24,000.$ Hence $y(n_v, \delta^{+})\ge y(\delta^{+},
\delta^{+})>0$ and
$$
{1.7n_v\ln\delta^+/\delta^+ \over 1+ \sqrt{6.8\ln\delta^+/\delta^+}}
> \ln(3.14n_v^{1.5}).
$$
We obtain
$$
2{\mathbf E}[|Q|]\le 2\sum_{v\in V(G)}{0.5n_v+1 \over
3.14n_v^{1.5}} \le {(\delta +3)n \over 3.14(\delta+1)^{1.5}}\le
{0.32 \over \sqrt{\delta+1}}n.
$$
Thus, (\ref{1}) yields
\begin{eqnarray*}
{\mathbf E}[g(V(G))]&\le& 2pn - n + {0.32n \over \sqrt{\delta +1}}\\
&=& {\sqrt{6.8\ln(\delta+1)}+ 0.32 \over \sqrt{\delta +1}}n,
\end{eqnarray*}
as required. The proof is complete. \qed

\begin{cor} \label{c2}
If $1,000 \le \delta \le 24,000,$ then
$$
\gamma_s(G) \le {\sqrt{\ln(\delta+1)(11.8-0.48\ln\delta)}+0.25 \over \sqrt{\delta+1}}n.
$$
\end{cor}

\pf It is similar to the proof of Corollary \ref{c1} if we put
$p=0.5 + \sqrt{(2.95-0.12\ln\delta)\ln{\delta^+}/\delta^+}$ and
consider the the following function for $1,000 \le \delta \le
24,000$:
$$
y(n_v, \delta^+)= {(2.95-0.12\ln\delta)n_v\ln\delta^+/\delta^+
\over 1+ \sqrt{(11.8-0.48\ln\delta)\ln\delta^+/\delta^+}} -
\ln(4.01n_v^{1.5}).
$$
\qed

\begin{cor} \label{c3}
If $230 \le \delta \le 1,000,$ then
$$
\gamma_s(G) \le {\sqrt{\ln(\delta+1)(18.16-1.4\ln\delta)}+0.25 \over \sqrt{\delta+1}}n.
$$
\end{cor}

\pf It is similar to the proof of Corollary \ref{c1} if we put
$p=0.5 + \sqrt{(4.54-0.35\ln\delta)\ln{\delta^+}/\delta^+}$ and
consider the following function for $230 \le \delta \le 1,000$:
$$
y(n_v, \delta^+)= {(4.54-0.35\ln\delta)n_v\ln\delta^+/\delta^+ \over
1+ \sqrt{(18.16-1.4\ln\delta)\ln\delta^+/\delta^+}} -
\ln(4.04n_v^{1.5}).
$$
\qed

We believe that F\"{u}redi--Mubayi's conjecture, saying that
$\ga_s(G) \le {Cn \over \sqrt\delta}$, is true for some small
constant $C$. However, as the Peterson graph shows, $C>1$, i.e.
the behaviour of the conjecture is not good for relatively small
values of $\delta$. Therefore, we propose the following rectified
conjecture, which, roughly speaking, consists of two functions for
`small' and `large' values of $\delta$.

\begin{conj}
For some $C\le 10$ and $\alpha$, $0.18\le \alpha < 0.21$,
$$
\ga_s(G) \le \min\Big\{{n\over \delta^\alpha}, {Cn\over
\sqrt\delta}\Big\}.
$$
\end{conj}

The above results imply that if $C=10$ and $\alpha=0.13$, then
this upper bound is true for all graphs with $\delta\le 96\times
10^4$.

\section{A Lower Bound for the Signed Domination Number}

The following theorem provides a lower bound for the signed
domination number of a graph $G$ depending on its order and a
parameter $\lambda$, which is  determined on the basis of the degree
sequence of $G$ (note that $\lambda$ may be equal to $0$, in this
case we put $\sum_{i=1}^{\lambda}=0$). This result improves the
bound of Theorem \ref{dunb_low} and, in some cases, it provides a
much better lower bound. For example, let us consider a graph $G$
consisting of two vertices of degree 5 and $n-2$ vertices of degree
3. Then, by Theorem \ref{dunb_low},
$$\gamma_s(G) \ge 0.25n-1,$$
while Theorem \ref{LowerBound} yields
$$\gamma_s(G) \ge 0.5n-1.$$

\begin{thm} \label{LowerBound}
Let $G$ be a graph with $n$ vertices and degrees $d_1\le d_2\le
... \le d_n$. Then
$$\gamma_s(G) \ge n-2\lambda,$$
where $\lambda\ge 0$ is the largest integer such that
$$
\sum_{i=1}^\lambda \left\lceil{d_i\over2}+1\right\rceil \le
\sum_{i=\lambda+1}^n \left\lfloor {d_i\over2}\right\rfloor.
$$
\end{thm}

\pf Let $f$ be a signed domination function of minimum weight of
the graph $G$. Let us denote
$$
X = \{v\in V(G): f(v)=-1\},
$$
and
$$
Y = \{v\in V(G): f(v)=1\}.
$$
We have
$$
\gamma_s(G) = f(V(G)) = |Y|-|X| = n - 2|X|.
$$
By definition, for any vertex $v$,
$$
f[v]=\sum_{u\in N[v]}f(u) \ge 1.
$$
Therefore, for all $v\in V(G)$,
$$
|N[v]\cap Y|-|N[v]\cap X| \ge 1.
$$
Using this inequality, we obtain
$$
|N[v]|=\deg (v)+1 = |N[v]\cap Y|+|N[v]\cap X| \le 2|N[v]\cap Y|-1.
$$
Hence
$$
|N[v]\cap Y| \ge {\deg(v)\over2} +1.
$$
Since $|N[v]\cap Y|$ is an integer, we conclude
$$
|N[v]\cap Y| \ge \left\lceil{\deg(v)\over2}\right\rceil +1
$$
and
$$
|N[v]\cap X| = \deg(v)+1 - |N[v]\cap Y| \le
\left\lfloor{\deg(v)\over2}\right\rfloor.
$$
Denote by $e_{X,Y}$ the number of edges between the parts $X$ and
$Y$. We have
$$
e_{X,Y} = \sum_{v\in X} |N[v]\cap Y| \ge \sum_{v\in X}
\left(\left\lceil{\deg(v)\over2}\right\rceil +1\right) \ge
\sum_{i=1}^{|X|} \left(\left\lceil{d_i\over2}\right\rceil
+1\right).
$$
Note that if $X=\emptyset$, then we put $\sum_{i=1}^0g(i)=0$. On
the other hand,
$$
e_{X,Y} = \sum_{v\in Y} |N[v]\cap X| \le \sum_{v\in Y}
\left\lfloor{\deg(v)\over2}\right\rfloor \le \sum_{i=n-|Y|+1}^n
\left\lfloor d_i/2\right\rfloor = \sum_{i=|X|+1}^n \lfloor
d_i/2\rfloor.
$$
Therefore, the following inequality holds:
$$
\sum_{i=1}^{|X|} \left(\left\lceil{d_i\over2}\right\rceil
+1\right) \le \sum_{i=|X|+1}^n
\left\lfloor{d_i\over2}\right\rfloor.
$$
Since $\lambda\ge 0$ is the largest integer such that
$$
\sum_{i=1}^{\lambda} \left(\left\lceil{d_i\over2}\right\rceil
+1\right) \le \sum_{i=\lambda+1}^n
\left\lfloor{d_i\over2}\right\rfloor,
$$
we conclude that
$$
|X| \le \lambda.
$$
Thus, $$\gamma_s(G) = n-2|X| \ge n-2\lambda.$$ The proof is
complete. \qed

Theorem \ref{LowerBound} immediately implies the following known
results:

\begin{cor} [\cite{haa1} and \cite{zha1}] \label{haas-zhang}
For any graph $G$,
$$
\gamma_s(G) \ge \left({\lceil0.5\delta\rceil -
\lfloor0.5\Delta\rfloor +1 \over \lceil0.5\delta\rceil +
\lfloor0.5\Delta\rfloor +1} \right)n.
$$
\end{cor}

Note that Haas and Wexler \cite{haa1} formulated the above bound
only for graphs with $\delta\ge 2$, while Zhang et al. \cite{zha1}
proved a weaker version without the ceiling and floor functions.

\begin{cor} [\cite{hen2}] \label{lb-Cor1}
If $\delta$ is odd and $G$ is $\delta$-regular, then
$$
\gamma_s(G) \ge {2n \over {\delta+1}}.
$$
\end{cor}

\begin{cor} [\cite{dun1}]  \label{lb-Cor2}
If $\delta$ is even and $G$ is $\delta$-regular, then
$$
\gamma_s(G) \ge {n \over {\delta+1}}.
$$
\end{cor}

Disjoint unions of complete graphs show that these lower bounds are
sharp whenever $n/(\delta+1)$ is an integer, and therefore the bound
of Theorem \ref{LowerBound} is sharp for regular graphs.


\end{document}